\nonstopmode
\input amstex
\input amsppt.sty   
%
\catcode`\@=11

\def\East#1#2{\setboxz@h{$\m@th\ssize\;{#1}\;\;$}%
 \setbox@ne\hbox{$\m@th\ssize\;{#2}\;\;$}\setbox\tw@\hbox{$\m@th#2$}%
 \dimen@\minaw@
 \ifdim\wdz@>\dimen@ \dimen@\wdz@ \fi  \ifdim\wd@ne>\dimen@ \dimen@\wd@ne \fi
 \ifdim\wd\tw@>\z@
  \mathrel{\mathop{\hbox to\dimen@{\rightarrowfill}}\limits^{#1}_{#2}}%
 \else
  \mathrel{\mathop{\hbox to\dimen@{\rightarrowfill}}\limits^{#1}}%
 \fi}
\def\West#1#2{\setboxz@h{$\m@th\ssize\;\;{#1}\;$}%
 \setbox@ne\hbox{$\m@th\ssize\;\;{#2}\;$}\setbox\tw@\hbox{$\m@th#2$}%
 \dimen@\minaw@
 \ifdim\wdz@>\dimen@ \dimen@\wdz@ \fi \ifdim\wd@ne>\dimen@ \dimen@\wd@ne \fi
 \ifdim\wd\tw@>\z@
  \mathrel{\mathop{\hbox to\dimen@{\leftarrowfill}}\limits^{#1}_{#2}}%
 \else
  \mathrel{\mathop{\hbox to\dimen@{\leftarrowfill}}\limits^{#1}}%
 \fi}
\font\arrow@i=lams1
\font\arrow@ii=lams2
\font\arrow@iii=lams3
\font\arrow@iv=lams4
\font\arrow@v=lams5
\newbox\zer@
\newdimen\standardcgap
\standardcgap=40\p@
\newdimen\hunit
\hunit=\tw@\p@
\newdimen\standardrgap
\standardrgap=32\p@
\newdimen\vunit
\vunit=1.6\p@
\def\Cgaps#1{\RIfM@
  \standardcgap=#1\standardcgap\relax \hunit=#1\hunit\relax
 \else \nonmatherr@\Cgaps \fi}
\def\Rgaps#1{\RIfM@
  \standardrgap=#1\standardrgap\relax \vunit=#1\vunit\relax
 \else \nonmatherr@\Rgaps \fi}
\newdimen\getdim@
\def\getcgap@#1{\ifcase#1\or\getdim@\z@\else\getdim@\standardcgap\fi}
\def\getrgap@#1{\ifcase#1\getdim@\z@\else\getdim@\standardrgap\fi}
\def\cgaps#1{\RIfM@
 \cgaps@{#1}\edef\getcgap@##1{\i@=##1\relax\the\toks@}\toks@{}\else
 \nonmatherr@\cgaps\fi}
\def\rgaps#1{\RIfM@
 \rgaps@{#1}\edef\getrgap@##1{\i@=##1\relax\the\toks@}\toks@{}\else
 \nonmatherr@\rgaps\fi}
\def\Gaps@@{\gaps@@}
\def\cgaps@#1{\toks@{\ifcase\i@\or\getdim@=\z@}%
 \gaps@@\standardcgap#1;\gaps@@\gaps@@
 \edef\next@{\the\toks@\noexpand\else\noexpand\getdim@\noexpand\standardcgap
  \noexpand\fi}%
 \toks@=\expandafter{\next@}}
\def\rgaps@#1{\toks@{\ifcase\i@\getdim@=\z@}%
 \gaps@@\standardrgap#1;\gaps@@\gaps@@
 \edef\next@{\the\toks@\noexpand\else\noexpand\getdim@\noexpand\standardrgap
  \noexpand\fi}%
 \toks@=\expandafter{\next@}}
\def\gaps@@#1#2;#3{\mgaps@#1#2\mgaps@
 \edef\next@{\the\toks@\noexpand\or\noexpand\getdim@
  \noexpand#1\the\mgapstoks@@}%
 \global\toks@=\expandafter{\next@}%
 \DN@{#3}%
 \ifx\next@\Gaps@@\gdef\next@##1\gaps@@{}\else
  \gdef\next@{\gaps@@#1#3}\fi\next@}
\def\mgaps@#1{\let\mgapsnext@#1\FN@\mgaps@@}
\def\mgaps@@{\ifx\next\space@\DN@. {\FN@\mgaps@@}\else
 \DN@.{\FN@\mgaps@@@}\fi\next@.}
\def\mgaps@@@{\ifx\next\w\let\next@\mgaps@@@@\else
 \let\next@\mgaps@@@@@\fi\next@}
\newtoks\mgapstoks@@
\def\mgaps@@@@@#1\mgaps@{\getdim@\mgapsnext@\getdim@#1\getdim@
 \edef\next@{\noexpand\getdim@\the\getdim@}%
 \mgapstoks@@=\expandafter{\next@}}
\def\mgaps@@@@\w#1#2\mgaps@{\mgaps@@@@@#2\mgaps@
 \setbox\zer@\hbox{$\m@th\hskip15\p@\tsize@#1$}%
 \dimen@\wd\zer@
 \ifdim\dimen@>\getdim@ \getdim@\dimen@ \fi
 \edef\next@{\noexpand\getdim@\the\getdim@}%
 \mgapstoks@@=\expandafter{\next@}}
\def\changewidth#1#2{\setbox\zer@\hbox{$\m@th#2}%
 \hbox to\wd\zer@{\hss$\m@th#1$\hss}}
\atdef@({\FN@\ARROW@}
\def\ARROW@{\ifx\next)\let\next@\OPTIONS@\else
 \DN@{\csname\string @(\endcsname}\fi\next@}
\newif\ifoptions@
\def\OPTIONS@){\ifoptions@\let\next@\relax\else
 \DN@{\options@true\begingroup\optioncodes@}\fi\next@}
\newif\ifN@
\newif\ifE@
\newif\ifNESW@
\newif\ifH@
\newif\ifV@
\newif\ifHshort@
\expandafter\def\csname\string @(\endcsname #1,#2){%
 \ifoptions@\let\next@\endgroup\else\let\next@\relax\fi\next@
 \N@false\E@false\H@false\V@false\Hshort@false
 \ifnum#1>\z@\E@true\fi
 \ifnum#1=\z@\V@true\tX@false\tY@false\a@false\fi
 \ifnum#2>\z@\N@true\fi
 \ifnum#2=\z@\H@true\tX@false\tY@false\a@false\ifshort@\Hshort@true\fi\fi
 \NESW@false
 \ifN@\ifE@\NESW@true\fi\else\ifE@\else\NESW@true\fi\fi
 \arrow@{#1}{#2}%
 \global\options@false
 \global\scount@\z@\global\tcount@\z@\global\arrcount@\z@
 \global\s@false\global\sxdimen@\z@\global\sydimen@\z@
 \global\tX@false\global\tXdimen@i\z@\global\tXdimen@ii\z@
 \global\tY@false\global\tYdimen@i\z@\global\tYdimen@ii\z@
 \global\a@false\global\exacount@\z@
 \global\x@false\global\xdimen@\z@
 \global\X@false\global\Xdimen@\z@
 \global\y@false\global\ydimen@\z@
 \global\Y@false\global\Ydimen@\z@
 \global\p@false\global\pdimen@\z@
 \global\label@ifalse\global\label@iifalse
 \global\dl@ifalse\global\ldimen@i\z@
 \global\dl@iifalse\global\ldimen@ii\z@
 \global\short@false\global\unshort@false}
\newif\iflabel@i
\newif\iflabel@ii
\newcount\scount@
\newcount\tcount@
\newcount\arrcount@
\newif\ifs@
\newdimen\sxdimen@
\newdimen\sydimen@
\newif\iftX@
\newdimen\tXdimen@i
\newdimen\tXdimen@ii
\newif\iftY@
\newdimen\tYdimen@i
\newdimen\tYdimen@ii
\newif\ifa@
\newcount\exacount@
\newif\ifx@
\newdimen\xdimen@
\newif\ifX@
\newdimen\Xdimen@
\newif\ify@
\newdimen\ydimen@
\newif\ifY@
\newdimen\Ydimen@
\newif\ifp@
\newdimen\pdimen@
\newif\ifdl@i
\newif\ifdl@ii
\newdimen\ldimen@i
\newdimen\ldimen@ii
\newif\ifshort@
\newif\ifunshort@
\def\zero@#1{\ifnum\scount@=\z@
 \if#1e\global\scount@\m@ne\else
 \if#1t\global\scount@\tw@\else
 \if#1h\global\scount@\thr@@\else
 \if#1'\global\scount@6 \else
 \if#1`\global\scount@7 \else
 \if#1(\global\scount@8 \else
 \if#1)\global\scount@9 \else
 \if#1s\global\scount@12 \else
 \if#1H\global\scount@13 \else
 \Err@{\Invalid@@ option \string\0}\fi\fi\fi\fi\fi\fi\fi\fi\fi
 \fi}
\def\one@#1{\ifnum\tcount@=\z@
 \if#1e\global\tcount@\m@ne\else
 \if#1h\global\tcount@\tw@\else
 \if#1t\global\tcount@\thr@@\else
 \if#1'\global\tcount@4 \else
 \if#1`\global\tcount@5 \else
 \if#1(\global\tcount@10 \else
 \if#1)\global\tcount@11 \else
 \if#1s\global\tcount@12 \else
 \if#1H\global\tcount@13 \else
 \Err@{\Invalid@@ option \string\1}\fi\fi\fi\fi\fi\fi\fi\fi\fi
 \fi}
\def\a@#1{\ifnum\arrcount@=\z@
 \if#10\global\arrcount@\m@ne\else
 \if#1+\global\arrcount@\@ne\else
 \if#1-\global\arrcount@\tw@\else
 \if#1=\global\arrcount@\thr@@\else
 \Err@{\Invalid@@ option \string\a}\fi\fi\fi\fi
 \fi}
\def\ds@(#1;#2){\ifs@\else
 \global\s@true
 \sxdimen@\hunit \global\sxdimen@#1\sxdimen@\relax
 \sydimen@\vunit \global\sydimen@#2\sydimen@\relax
 \fi}
\def\dtX@(#1;#2){\iftX@\else
 \global\tX@true
 \tXdimen@i\hunit \global\tXdimen@i#1\tXdimen@i\relax
 \tXdimen@ii\vunit \global\tXdimen@ii#2\tXdimen@ii\relax
 \fi}
\def\dtY@(#1;#2){\iftY@\else
 \global\tY@true
 \tYdimen@i\hunit \global\tYdimen@i#1\tYdimen@i\relax
 \tYdimen@ii\vunit \global\tYdimen@ii#2\tYdimen@ii\relax
 \fi}
\def\da@#1{\ifa@\else\global\a@true\global\exacount@#1\relax\fi}
\def\dx@#1{\ifx@\else
 \global\x@true
 \xdimen@\hunit \global\xdimen@#1\xdimen@\relax
 \fi}
\def\dX@#1{\ifX@\else
 \global\X@true
 \Xdimen@\hunit \global\Xdimen@#1\Xdimen@\relax
 \fi}
\def\dy@#1{\ify@\else
 \global\y@true
 \ydimen@\vunit \global\ydimen@#1\ydimen@\relax
 \fi}
\def\dY@#1{\ifY@\else
 \global\Y@true
 \Ydimen@\vunit \global\Ydimen@#1\Ydimen@\relax
 \fi}
\def\p@@#1{\ifp@\else
 \global\p@true
 \pdimen@\hunit \divide\pdimen@\tw@ \global\pdimen@#1\pdimen@\relax
 \fi}
\def\L@#1{\iflabel@i\else
 \global\label@itrue \gdef\label@i{#1}%
 \fi}
\def\l@#1{\iflabel@ii\else
 \global\label@iitrue \gdef\label@ii{#1}%
 \fi}
\def\dL@#1{\ifdl@i\else
 \global\dl@itrue \ldimen@i\hunit \global\ldimen@i#1\ldimen@i\relax
 \fi}
\def\dl@#1{\ifdl@ii\else
 \global\dl@iitrue \ldimen@ii\hunit \global\ldimen@ii#1\ldimen@ii\relax
 \fi}
\def\s@{\ifunshort@\else\global\short@true\fi}
\def\uns@{\ifshort@\else\global\unshort@true\global\short@false\fi}
\def\optioncodes@{\let\0\zero@\let\1\one@\let\a\a@\let\ds\ds@\let\dtX\dtX@
 \let\dtY\dtY@\let\da\da@\let\dx\dx@\let\dX\dX@\let\dY\dY@\let\dy\dy@
 \let\p\p@@\let\L\L@\let\l\l@\let\dL\dL@\let\dl\dl@\let\s\s@\let\uns\uns@}
\def\slopes@{\\161\\152\\143\\134\\255\\126\\357\\238\\349\\45{10}\\56{11}%
 \\11{12}\\65{13}\\54{14}\\43{15}\\32{16}\\53{17}\\21{18}\\52{19}\\31{20}%
 \\41{21}\\51{22}\\61{23}}
\newcount\tan@i
\newcount\tan@ip
\newcount\tan@ii
\newcount\tan@iip
\newdimen\slope@i
\newdimen\slope@ip
\newdimen\slope@ii
\newdimen\slope@iip
\newcount\angcount@
\newcount\extracount@
\def\slope@{{\slope@i=\secondy@ \advance\slope@i-\firsty@
 \ifN@\else\multiply\slope@i\m@ne\fi
 \slope@ii=\secondx@ \advance\slope@ii-\firstx@
 \ifE@\else\multiply\slope@ii\m@ne\fi
 \ifdim\slope@ii<\z@
  \global\tan@i6 \global\tan@ii\@ne \global\angcount@23
 \else
  \dimen@\slope@i \multiply\dimen@6
  \ifdim\dimen@<\slope@ii
   \global\tan@i\@ne \global\tan@ii6 \global\angcount@\@ne
  \else
   \dimen@\slope@ii \multiply\dimen@6
   \ifdim\dimen@<\slope@i
    \global\tan@i6 \global\tan@ii\@ne \global\angcount@23
   \else
    \tan@ip\z@ \tan@iip \@ne
    \def\\##1##2##3{\global\angcount@=##3\relax
     \slope@ip\slope@i \slope@iip\slope@ii
     \multiply\slope@iip##1\relax \multiply\slope@ip##2\relax
     \ifdim\slope@iip<\slope@ip
      \global\tan@ip=##1\relax \global\tan@iip=##2\relax
     \else
      \global\tan@i=##1\relax \global\tan@ii=##2\relax
      \def\\####1####2####3{}%
     \fi}%
    \slopes@
    \slope@i=\secondy@ \advance\slope@i-\firsty@
    \ifN@\else\multiply\slope@i\m@ne\fi
    \multiply\slope@i\tan@ii \multiply\slope@i\tan@iip \multiply\slope@i\tw@
    \count@\tan@i \multiply\count@\tan@iip
    \extracount@\tan@ip \multiply\extracount@\tan@ii
    \advance\count@\extracount@
    \slope@ii=\secondx@ \advance\slope@ii-\firstx@
    \ifE@\else\multiply\slope@ii\m@ne\fi
    \multiply\slope@ii\count@
    \ifdim\slope@i<\slope@ii
     \global\tan@i=\tan@ip \global\tan@ii=\tan@iip
     \global\advance\angcount@\m@ne
    \fi
   \fi
  \fi
 \fi}%
}
\def\slope@a#1{{\def\\##1##2##3{\ifnum##3=#1\global\tan@i=##1\relax
 \global\tan@ii=##2\relax\fi}\slopes@}}
\newcount\i@
\newcount\j@
\newcount\colcount@
\newcount\Colcount@
\newcount\tcolcount@
\newdimen\rowht@
\newdimen\rowdp@
\newcount\rowcount@
\newcount\Rowcount@
\newcount\maxcolrow@
\newtoks\colwidthtoks@
\newtoks\Rowheighttoks@
\newtoks\Rowdepthtoks@
\newtoks\widthtoks@
\newtoks\Widthtoks@
\newtoks\heighttoks@
\newtoks\Heighttoks@
\newtoks\depthtoks@
\newtoks\Depthtoks@
\newif\iffirstnewCDcr@
\def\dotoks@i{%
 \global\widthtoks@=\expandafter{\the\widthtoks@\else\getdim@\z@\fi}%
 \global\heighttoks@=\expandafter{\the\heighttoks@\else\getdim@\z@\fi}%
 \global\depthtoks@=\expandafter{\the\depthtoks@\else\getdim@\z@\fi}}
\def\dotoks@ii{%
 \global\widthtoks@{\ifcase\j@}%
 \global\heighttoks@{\ifcase\j@}%
 \global\depthtoks@{\ifcase\j@}}
\def\prenewCD@#1\endnewCD{\setbox\zer@
 \vbox{%
  \def\arrow@##1##2{{}}%
  \rowcount@\m@ne \colcount@\z@ \Colcount@\z@
  \firstnewCDcr@true \toks@{}%
  \widthtoks@{\ifcase\j@}%
  \Widthtoks@{\ifcase\i@}%
  \heighttoks@{\ifcase\j@}%
  \Heighttoks@{\ifcase\i@}%
  \depthtoks@{\ifcase\j@}%
  \Depthtoks@{\ifcase\i@}%
  \Rowheighttoks@{\ifcase\i@}%
  \Rowdepthtoks@{\ifcase\i@}%
  \Let@
  \everycr{%
   \noalign{%
    \global\advance\rowcount@\@ne
    \ifnum\colcount@<\Colcount@
    \else
     \global\Colcount@=\colcount@ \global\maxcolrow@=\rowcount@
    \fi
    \global\colcount@\z@
    \iffirstnewCDcr@
     \global\firstnewCDcr@false
    \else
     \edef\next@{\the\Rowheighttoks@\noexpand\or\noexpand\getdim@\the\rowht@}%
      \global\Rowheighttoks@=\expandafter{\next@}%
     \edef\next@{\the\Rowdepthtoks@\noexpand\or\noexpand\getdim@\the\rowdp@}%
      \global\Rowdepthtoks@=\expandafter{\next@}%
     \global\rowht@\z@ \global\rowdp@\z@
     \dotoks@i
     \edef\next@{\the\Widthtoks@\noexpand\or\the\widthtoks@}%
      \global\Widthtoks@=\expandafter{\next@}%
     \edef\next@{\the\Heighttoks@\noexpand\or\the\heighttoks@}%
      \global\Heighttoks@=\expandafter{\next@}%
     \edef\next@{\the\Depthtoks@\noexpand\or\the\depthtoks@}%
      \global\Depthtoks@=\expandafter{\next@}%
     \dotoks@ii
    \fi}}%
  \tabskip\z@
  \halign{&\setbox\zer@\hbox{\vrule height10\p@ width\z@ depth\z@
   $\m@th\displaystyle{##}$}\copy\zer@
   \ifdim\ht\zer@>\rowht@ \global\rowht@\ht\zer@ \fi
   \ifdim\dp\zer@>\rowdp@ \global\rowdp@\dp\zer@ \fi
   \global\advance\colcount@\@ne
   \edef\next@{\the\widthtoks@\noexpand\or\noexpand\getdim@\the\wd\zer@}%
    \global\widthtoks@=\expandafter{\next@}%
   \edef\next@{\the\heighttoks@\noexpand\or\noexpand\getdim@\the\ht\zer@}%
    \global\heighttoks@=\expandafter{\next@}%
   \edef\next@{\the\depthtoks@\noexpand\or\noexpand\getdim@\the\dp\zer@}%
    \global\depthtoks@=\expandafter{\next@}%
   \cr#1\crcr}}%
 \Rowcount@=\rowcount@
 \global\Widthtoks@=\expandafter{\the\Widthtoks@\fi\relax}%
 \edef\Width@##1##2{\i@=##1\relax\j@=##2\relax\the\Widthtoks@}%
 \global\Heighttoks@=\expandafter{\the\Heighttoks@\fi\relax}%
 \edef\Height@##1##2{\i@=##1\relax\j@=##2\relax\the\Heighttoks@}%
 \global\Depthtoks@=\expandafter{\the\Depthtoks@\fi\relax}%
 \edef\Depth@##1##2{\i@=##1\relax\j@=##2\relax\the\Depthtoks@}%
 \edef\next@{\the\Rowheighttoks@\noexpand\fi\relax}%
 \global\Rowheighttoks@=\expandafter{\next@}%
 \edef\Rowheight@##1{\i@=##1\relax\the\Rowheighttoks@}%
 \edef\next@{\the\Rowdepthtoks@\noexpand\fi\relax}%
 \global\Rowdepthtoks@=\expandafter{\next@}%
 \edef\Rowdepth@##1{\i@=##1\relax\the\Rowdepthtoks@}%
 \colwidthtoks@{\fi}%
 \setbox\zer@\vbox{%
  \unvbox\zer@
  \count@\rowcount@
  \loop
   \unskip\unpenalty
   \setbox\zer@\lastbox
   \ifnum\count@>\maxcolrow@ \advance\count@\m@ne
   \repeat
  \hbox{%
   \unhbox\zer@
   \count@\z@
   \loop
    \unskip
    \setbox\zer@\lastbox
    \edef\next@{\noexpand\or\noexpand\getdim@\the\wd\zer@\the\colwidthtoks@}%
     \global\colwidthtoks@=\expandafter{\next@}%
    \advance\count@\@ne
    \ifnum\count@<\Colcount@
    \repeat}}%
 \edef\next@{\noexpand\ifcase\noexpand\i@\the\colwidthtoks@}%
  \global\colwidthtoks@=\expandafter{\next@}%
 \edef\Colwidth@##1{\i@=##1\relax\the\colwidthtoks@}%
 \colwidthtoks@{}\Rowheighttoks@{}\Rowdepthtoks@{}\widthtoks@{}%
 \Widthtoks@{}\heighttoks@{}\Heighttoks@{}\depthtoks@{}\Depthtoks@{}%
}
\newcount\xoff@
\newcount\yoff@
\newcount\endcount@
\newcount\rcount@
\newdimen\firstx@
\newdimen\firsty@
\newdimen\secondx@
\newdimen\secondy@
\newdimen\tocenter@
\newdimen\charht@
\newdimen\charwd@
\def\outside@{\Err@{This arrow points outside the \string\newCD}}
\newif\ifsvertex@
\newif\iftvertex@
\def\arrow@#1#2{\xoff@=#1\relax\yoff@=#2\relax
 \count@\rowcount@ \advance\count@-\yoff@
 \ifnum\count@<\@ne \outside@ \else \ifnum\count@>\Rowcount@ \outside@ \fi\fi
 \count@\colcount@ \advance\count@\xoff@
 \ifnum\count@<\@ne \outside@ \else \ifnum\count@>\Colcount@ \outside@\fi\fi
 \tcolcount@\colcount@ \advance\tcolcount@\xoff@
 \Width@\rowcount@\colcount@ \tocenter@=-\getdim@ \divide\tocenter@\tw@
 \ifdim\getdim@=\z@
  \firstx@\z@ \firsty@\mathaxis@ \svertex@true
 \else
  \svertex@false
  \ifHshort@
   \Colwidth@\colcount@
    \ifE@ \firstx@=.5\getdim@ \else \firstx@=-.5\getdim@ \fi
  \else
   \ifE@ \firstx@=\getdim@ \else \firstx@=-\getdim@ \fi
   \divide\firstx@\tw@
  \fi
  \ifE@
   \ifH@ \advance\firstx@\thr@@\p@ \else \advance\firstx@-\thr@@\p@ \fi
  \else
   \ifH@ \advance\firstx@-\thr@@\p@ \else \advance\firstx@\thr@@\p@ \fi
  \fi
  \ifN@
   \Height@\rowcount@\colcount@ \firsty@=\getdim@
   \ifV@ \advance\firsty@\thr@@\p@ \fi
  \else
   \ifV@
    \Depth@\rowcount@\colcount@ \firsty@=-\getdim@
    \advance\firsty@-\thr@@\p@
   \else
    \firsty@\z@
   \fi
  \fi
 \fi
 \ifV@
 \else
  \Colwidth@\colcount@
  \ifE@ \secondx@=\getdim@ \else \secondx@=-\getdim@ \fi
  \divide\secondx@\tw@
  \ifE@ \else \getcgap@\colcount@ \advance\secondx@-\getdim@ \fi
  \endcount@=\colcount@ \advance\endcount@\xoff@
  \count@=\colcount@
  \ifE@
   \advance\count@\@ne
   \loop
    \ifnum\count@<\endcount@
    \Colwidth@\count@ \advance\secondx@\getdim@
    \getcgap@\count@ \advance\secondx@\getdim@
    \advance\count@\@ne
    \repeat
  \else
   \advance\count@\m@ne
   \loop
    \ifnum\count@>\endcount@
    \Colwidth@\count@ \advance\secondx@-\getdim@
    \getcgap@\count@ \advance\secondx@-\getdim@
    \advance\count@\m@ne
    \repeat
  \fi
  \Colwidth@\count@ \divide\getdim@\tw@
  \ifHshort@
  \else
   \ifE@ \advance\secondx@\getdim@ \else \advance\secondx@-\getdim@ \fi
  \fi
  \ifE@ \getcgap@\count@ \advance\secondx@\getdim@ \fi
  \rcount@\rowcount@ \advance\rcount@-\yoff@
  \Width@\rcount@\count@ \divide\getdim@\tw@
  \tvertex@false
  \ifH@\ifdim\getdim@=\z@\tvertex@true\Hshort@false\fi\fi
  \ifHshort@
  \else
   \ifE@ \advance\secondx@-\getdim@ \else \advance\secondx@\getdim@ \fi
  \fi
  \iftvertex@
   \advance\secondx@.4\p@
  \else
   \ifE@ \advance\secondx@-\thr@@\p@ \else \advance\secondx@\thr@@\p@ \fi
  \fi
 \fi
 \ifH@
 \else
  \ifN@
   \Rowheight@\rowcount@ \secondy@\getdim@
  \else
   \Rowdepth@\rowcount@ \secondy@-\getdim@
   \getrgap@\rowcount@ \advance\secondy@-\getdim@
  \fi
  \endcount@=\rowcount@ \advance\endcount@-\yoff@
  \count@=\rowcount@
  \ifN@
   \advance\count@\m@ne
   \loop
    \ifnum\count@>\endcount@
    \Rowheight@\count@ \advance\secondy@\getdim@
    \Rowdepth@\count@ \advance\secondy@\getdim@
    \getrgap@\count@ \advance\secondy@\getdim@
    \advance\count@\m@ne
    \repeat
  \else
   \advance\count@\@ne
   \loop
    \ifnum\count@<\endcount@
    \Rowheight@\count@ \advance\secondy@-\getdim@
    \Rowdepth@\count@ \advance\secondy@-\getdim@
    \getrgap@\count@ \advance\secondy@-\getdim@
    \advance\count@\@ne
    \repeat
  \fi
  \tvertex@false
  \ifV@\Width@\count@\colcount@\ifdim\getdim@=\z@\tvertex@true\fi\fi
  \ifN@
   \getrgap@\count@ \advance\secondy@\getdim@
   \Rowdepth@\count@ \advance\secondy@\getdim@
   \iftvertex@
    \advance\secondy@\mathaxis@
   \else
    \Depth@\count@\tcolcount@ \advance\secondy@-\getdim@
    \advance\secondy@-\thr@@\p@
   \fi
  \else
   \Rowheight@\count@ \advance\secondy@-\getdim@
   \iftvertex@
    \advance\secondy@\mathaxis@
   \else
    \Height@\count@\tcolcount@ \advance\secondy@\getdim@
    \advance\secondy@\thr@@\p@
   \fi
  \fi
 \fi
 \ifV@\else\advance\firstx@\sxdimen@\fi
 \ifH@\else\advance\firsty@\sydimen@\fi
 \iftX@
  \advance\secondy@\tXdimen@ii
  \advance\secondx@\tXdimen@i
  \slope@
 \else
  \iftY@
   \advance\secondy@\tYdimen@ii
   \advance\secondx@\tYdimen@i
   \slope@
   \secondy@=\secondx@ \advance\secondy@-\firstx@
   \ifNESW@ \else \multiply\secondy@\m@ne \fi
   \multiply\secondy@\tan@i \divide\secondy@\tan@ii \advance\secondy@\firsty@
  \else
   \ifa@
    \slope@
    \ifNESW@ \global\advance\angcount@\exacount@ \else
      \global\advance\angcount@-\exacount@ \fi
    \ifnum\angcount@>23 \angcount@23 \fi
    \ifnum\angcount@<\@ne \angcount@\@ne \fi
    \slope@a\angcount@
    \ifY@
     \advance\secondy@\Ydimen@
    \else
     \ifX@
      \advance\secondx@\Xdimen@
      \dimen@\secondx@ \advance\dimen@-\firstx@
      \ifNESW@\else\multiply\dimen@\m@ne\fi
      \multiply\dimen@\tan@i \divide\dimen@\tan@ii
      \advance\dimen@\firsty@ \secondy@=\dimen@
     \fi
    \fi
   \else
    \ifH@\else\ifV@\else\slope@\fi\fi
   \fi
  \fi
 \fi
 \ifH@\else\ifV@\else\ifsvertex@\else
  \dimen@=6\p@ \multiply\dimen@\tan@ii
  \count@=\tan@i \advance\count@\tan@ii \divide\dimen@\count@
  \ifE@ \advance\firstx@\dimen@ \else \advance\firstx@-\dimen@ \fi
  \multiply\dimen@\tan@i \divide\dimen@\tan@ii
  \ifN@ \advance\firsty@\dimen@ \else \advance\firsty@-\dimen@ \fi
 \fi\fi\fi
 \ifp@
  \ifH@\else\ifV@\else
   \getcos@\pdimen@ \advance\firsty@\dimen@ \advance\secondy@\dimen@
   \ifNESW@ \advance\firstx@-\dimen@ii \else \advance\firstx@\dimen@ii \fi
  \fi\fi
 \fi
 \ifH@\else\ifV@\else
  \ifnum\tan@i>\tan@ii
   \charht@=10\p@ \charwd@=10\p@
   \multiply\charwd@\tan@ii \divide\charwd@\tan@i
  \else
   \charwd@=10\p@ \charht@=10\p@
   \divide\charht@\tan@ii \multiply\charht@\tan@i
  \fi
  \ifnum\tcount@=\thr@@
   \ifN@ \advance\secondy@-.3\charht@ \else\advance\secondy@.3\charht@ \fi
  \fi
  \ifnum\scount@=\tw@
   \ifE@ \advance\firstx@.3\charht@ \else \advance\firstx@-.3\charht@ \fi
  \fi
  \ifnum\tcount@=12
   \ifN@ \advance\secondy@-\charht@ \else \advance\secondy@\charht@ \fi
  \fi
  \iftY@
  \else
   \ifa@
    \ifX@
    \else
     \secondx@\secondy@ \advance\secondx@-\firsty@
     \ifNESW@\else\multiply\secondx@\m@ne\fi
     \multiply\secondx@\tan@ii \divide\secondx@\tan@i
     \advance\secondx@\firstx@
    \fi
   \fi
  \fi
 \fi\fi
 \ifH@\harrow@\else\ifV@\varrow@\else\arrow@@\fi\fi}
\newdimen\mathaxis@
\mathaxis@90\p@ \divide\mathaxis@36
\def\harrow@b{\ifE@\hskip\tocenter@\hskip\firstx@\fi}
\def\harrow@bb{\ifE@\hskip\xdimen@\else\hskip\Xdimen@\fi}
\def\harrow@e{\ifE@\else\hskip-\firstx@\hskip-\tocenter@\fi}
\def\harrow@ee{\ifE@\hskip-\Xdimen@\else\hskip-\xdimen@\fi}
\def\harrow@{\dimen@\secondx@\advance\dimen@-\firstx@
 \ifE@ \let\next@\rlap \else  \multiply\dimen@\m@ne \let\next@\llap \fi
 \next@{%
  \harrow@b
  \smash{\raise\pdimen@\hbox to\dimen@
   {\harrow@bb\arrow@ii
    \ifnum\arrcount@=\m@ne \else \ifnum\arrcount@=\thr@@ \else
     \ifE@
      \ifnum\scount@=\m@ne
      \else
       \ifcase\scount@\or\or\char118 \or\char117 \or\or\or\char119 \or
       \char120 \or\char121 \or\char122 \or\or\or\arrow@i\char125 \or
       \char117 \hskip\thr@@\p@\char117 \hskip-\thr@@\p@\fi
      \fi
     \else
      \ifnum\tcount@=\m@ne
      \else
       \ifcase\tcount@\char117 \or\or\char117 \or\char118 \or\char119 \or
       \char120\or\or\or\or\or\char121 \or\char122 \or\arrow@i\char125
       \or\char117 \hskip\thr@@\p@\char117 \hskip-\thr@@\p@\fi
      \fi
     \fi
    \fi\fi
    \dimen@\mathaxis@ \advance\dimen@.2\p@
    \dimen@ii\mathaxis@ \advance\dimen@ii-.2\p@
    \ifnum\arrcount@=\m@ne
     \let\leads@\null
    \else
     \ifcase\arrcount@
      \def\leads@{\hrule height\dimen@ depth-\dimen@ii}\or
      \def\leads@{\hrule height\dimen@ depth-\dimen@ii}\or
      \def\leads@{\hbox to10\p@{%
       \leaders\hrule height\dimen@ depth-\dimen@ii\hfil
       \hfil
      \leaders\hrule height\dimen@ depth-\dimen@ii\hskip\z@ plus2fil\relax
       \hfil
       \leaders\hrule height\dimen@ depth-\dimen@ii\hfil}}\or
     \def\leads@{\hbox{\hbox to10\p@{\dimen@\mathaxis@ \advance\dimen@1.2\p@
       \dimen@ii\dimen@ \advance\dimen@ii-.4\p@
       \leaders\hrule height\dimen@ depth-\dimen@ii\hfil}%
       \kern-10\p@
       \hbox to10\p@{\dimen@\mathaxis@ \advance\dimen@-1.2\p@
       \dimen@ii\dimen@ \advance\dimen@ii-.4\p@
       \leaders\hrule height\dimen@ depth-\dimen@ii\hfil}}}\fi
    \fi
    \cleaders\leads@\hfil
    \ifnum\arrcount@=\m@ne\else\ifnum\arrcount@=\thr@@\else
     \arrow@i
     \ifE@
      \ifnum\tcount@=\m@ne
      \else
       \ifcase\tcount@\char119 \or\or\char119 \or\char120 \or\char121 \or
       \char122 \or \or\or\or\or\char123\or\char124 \or
       \char125 \or\char119 \hskip-\thr@@\p@\char119 \hskip\thr@@\p@\fi
      \fi
     \else
      \ifcase\scount@\or\or\char120 \or\char119 \or\or\or\char121 \or\char122
      \or\char123 \or\char124 \or\or\or\char125 \or
      \char119 \hskip-\thr@@\p@\char119 \hskip\thr@@\p@\fi
     \fi
    \fi\fi
    \harrow@ee}}%
  \harrow@e}%
 \iflabel@i
  \dimen@ii\z@ \setbox\zer@\hbox{$\m@th\tsize@@\label@i$}%
  \ifnum\arrcount@=\m@ne
  \else
   \advance\dimen@ii\mathaxis@
   \advance\dimen@ii\dp\zer@ \advance\dimen@ii\tw@\p@
   \ifnum\arrcount@=\thr@@ \advance\dimen@ii\tw@\p@ \fi
  \fi
  \advance\dimen@ii\pdimen@
  \next@{\harrow@b\smash{\raise\dimen@ii\hbox to\dimen@
   {\harrow@bb\hskip\tw@\ldimen@i\hfil\box\zer@\hfil\harrow@ee}}\harrow@e}%
 \fi
 \iflabel@ii
  \ifnum\arrcount@=\m@ne
  \else
   \setbox\zer@\hbox{$\m@th\tsize@\label@ii$}%
   \dimen@ii-\ht\zer@ \advance\dimen@ii-\tw@\p@
   \ifnum\arrcount@=\thr@@ \advance\dimen@ii-\tw@\p@ \fi
   \advance\dimen@ii\mathaxis@ \advance\dimen@ii\pdimen@
   \next@{\harrow@b\smash{\raise\dimen@ii\hbox to\dimen@
    {\harrow@bb\hskip\tw@\ldimen@ii\hfil\box\zer@\hfil\harrow@ee}}\harrow@e}%
  \fi
 \fi}
\let\tsize@\tsize
\def\tsizenewCDlabels{\let\tsize@\tsize}
\def\ssizenewCDlabels{\let\tsize@\ssize}
\def\tsize@@{\ifnum\arrcount@=\m@ne\else\tsize@\fi}
\def\varrow@{\dimen@\secondy@ \advance\dimen@-\firsty@
 \ifN@ \else \multiply\dimen@\m@ne \fi
 \setbox\zer@\vbox to\dimen@
  {\ifN@ \vskip-\Ydimen@ \else \vskip\ydimen@ \fi
   \ifnum\arrcount@=\m@ne\else\ifnum\arrcount@=\thr@@\else
    \hbox{\arrow@iii
     \ifN@
      \ifnum\tcount@=\m@ne
      \else
       \ifcase\tcount@\char117 \or\or\char117 \or\char118 \or\char119 \or
       \char120 \or\or\or\or\or\char121 \or\char122 \or\char123 \or
       \vbox{\hbox{\char117 }\nointerlineskip\vskip\thr@@\p@
       \hbox{\char117 }\vskip-\thr@@\p@}\fi
      \fi
     \else
      \ifcase\scount@\or\or\char118 \or\char117 \or\or\or\char119 \or
      \char120 \or\char121 \or\char122 \or\or\or\char123 \or
      \vbox{\hbox{\char117 }\nointerlineskip\vskip\thr@@\p@
      \hbox{\char117 }\vskip-\thr@@\p@}\fi
     \fi}%
    \nointerlineskip
   \fi\fi
   \ifnum\arrcount@=\m@ne
    \let\leads@\null
   \else
    \ifcase\arrcount@\let\leads@\vrule\or\let\leads@\vrule\or
    \def\leads@{\vbox to10\p@{%
     \hrule height 1.67\p@ depth\z@ width.4\p@
     \vfil
     \hrule height 3.33\p@ depth\z@ width.4\p@
     \vfil
     \hrule height 1.67\p@ depth\z@ width.4\p@}}\or
    \def\leads@{\hbox{\vrule height\p@\hskip\tw@\p@\vrule}}\fi
   \fi
  \cleaders\leads@\vfill\nointerlineskip
   \ifnum\arrcount@=\m@ne\else\ifnum\arrcount@=\thr@@\else
    \hbox{\arrow@iv
     \ifN@
      \ifcase\scount@\or\or\char118 \or\char117 \or\or\or\char119 \or
      \char120 \or\char121 \or\char122 \or\or\or\arrow@iii\char123 \or
      \vbox{\hbox{\char117 }\nointerlineskip\vskip-\thr@@\p@
      \hbox{\char117 }\vskip\thr@@\p@}\fi
     \else
      \ifnum\tcount@=\m@ne
      \else
       \ifcase\tcount@\char117 \or\or\char117 \or\char118 \or\char119 \or
       \char120 \or\or\or\or\or\char121 \or\char122 \or\arrow@iii\char123 \or
       \vbox{\hbox{\char117 }\nointerlineskip\vskip-\thr@@\p@
       \hbox{\char117 }\vskip\thr@@\p@}\fi
      \fi
     \fi}%
   \fi\fi
   \ifN@\vskip\ydimen@\else\vskip-\Ydimen@\fi}%
 \ifN@
  \dimen@ii\firsty@
 \else
  \dimen@ii-\firsty@ \advance\dimen@ii\ht\zer@ \multiply\dimen@ii\m@ne
 \fi
 \rlap{\smash{\hskip\tocenter@ \hskip\pdimen@ \raise\dimen@ii \box\zer@}}%
 \iflabel@i
  \setbox\zer@\vbox to\dimen@{\vfil
   \hbox{$\m@th\tsize@@\label@i$}\vskip\tw@\ldimen@i\vfil}%
  \rlap{\smash{\hskip\tocenter@ \hskip\pdimen@
  \ifnum\arrcount@=\m@ne \let\next@\relax \else \let\next@\llap \fi
  \next@{\raise\dimen@ii\hbox{\ifnum\arrcount@=\m@ne \hskip-.5\wd\zer@ \fi
   \box\zer@ \ifnum\arrcount@=\m@ne \else \hskip\tw@\p@ \fi}}}}%
 \fi
 \iflabel@ii
  \ifnum\arrcount@=\m@ne
  \else
   \setbox\zer@\vbox to\dimen@{\vfil
    \hbox{$\m@th\tsize@\label@ii$}\vskip\tw@\ldimen@ii\vfil}%
   \rlap{\smash{\hskip\tocenter@ \hskip\pdimen@
   \rlap{\raise\dimen@ii\hbox{\ifnum\arrcount@=\thr@@ \hskip4.5\p@ \else
    \hskip2.5\p@ \fi\box\zer@}}}}%
  \fi
 \fi
}
\newdimen\goal@
\newdimen\shifted@
\newcount\Tcount@
\newcount\Scount@
\newbox\shaft@
\newcount\slcount@
\def\getcos@#1{%
 \ifnum\tan@i<\tan@ii
  \dimen@#1%
  \ifnum\slcount@<8 \count@9 \else \ifnum\slcount@<12 \count@8 \else
   \count@7 \fi\fi
  \multiply\dimen@\count@ \divide\dimen@10
  \dimen@ii\dimen@ \multiply\dimen@ii\tan@i \divide\dimen@ii\tan@ii
 \else
  \dimen@ii#1%
  \count@-\slcount@ \advance\count@24
  \ifnum\count@<8 \count@9 \else \ifnum\count@<12 \count@8
   \else\count@7 \fi\fi
  \multiply\dimen@ii\count@ \divide\dimen@ii10
  \dimen@\dimen@ii \multiply\dimen@\tan@ii \divide\dimen@\tan@i
 \fi}
\newdimen\adjust@
\def\Nnext@{\ifN@\let\next@\raise\else\let\next@\lower\fi}
\def\arrow@@{\slcount@\angcount@
 \ifNESW@
  \ifnum\angcount@<10
   \let\arrowfont@=\arrow@i \advance\angcount@\m@ne \multiply\angcount@13
  \else
   \ifnum\angcount@<19
    \let\arrowfont@=\arrow@ii \advance\angcount@-10 \multiply\angcount@13
   \else
    \let\arrowfont@=\arrow@iii \advance\angcount@-19 \multiply\angcount@13
  \fi\fi
  \Tcount@\angcount@
 \else
  \ifnum\angcount@<5
   \let\arrowfont@=\arrow@iii \advance\angcount@\m@ne \multiply\angcount@13
   \advance\angcount@65
  \else
   \ifnum\angcount@<14
    \let\arrowfont@=\arrow@iv \advance\angcount@-5 \multiply\angcount@13
   \else
    \ifnum\angcount@<23
     \let\arrowfont@=\arrow@v \advance\angcount@-14 \multiply\angcount@13
    \else
     \let\arrowfont@=\arrow@i \angcount@=117
  \fi\fi\fi
  \ifnum\angcount@=117 \Tcount@=115 \else\Tcount@\angcount@ \fi
 \fi
 \Scount@\Tcount@
 \ifE@
  \ifnum\tcount@=\z@ \advance\Tcount@\tw@ \else\ifnum\tcount@=13
   \advance\Tcount@\tw@ \else \advance\Tcount@\tcount@ \fi\fi
  \ifnum\scount@=\z@ \else \ifnum\scount@=13 \advance\Scount@\thr@@ \else
   \advance\Scount@\scount@ \fi\fi
 \else
  \ifcase\tcount@\advance\Tcount@\thr@@\or\or\advance\Tcount@\thr@@\or
  \advance\Tcount@\tw@\or\advance\Tcount@6 \or\advance\Tcount@7
  \or\or\or\or\or \advance\Tcount@8 \or\advance\Tcount@9 \or
  \advance\Tcount@12 \or\advance\Tcount@\thr@@\fi
  \ifcase\scount@\or\or\advance\Scount@\thr@@\or\advance\Scount@\tw@\or
  \or\or\advance\Scount@4 \or\advance\Scount@5 \or\advance\Scount@10
  \or\advance\Scount@11 \or\or\or\advance\Scount@12 \or\advance
  \Scount@\tw@\fi
 \fi
 \ifcase\arrcount@\or\or\advance\angcount@\@ne\else\fi
 \ifN@ \shifted@=\firsty@ \else\shifted@=-\firsty@ \fi
 \ifE@ \else\advance\shifted@\charht@ \fi
 \goal@=\secondy@ \advance\goal@-\firsty@
 \ifN@\else\multiply\goal@\m@ne\fi
 \setbox\shaft@\hbox{\arrowfont@\char\angcount@}%
 \ifnum\arrcount@=\thr@@
  \getcos@{1.5\p@}%
  \setbox\shaft@\hbox to\wd\shaft@{\arrowfont@
   \rlap{\hskip\dimen@ii
    \smash{\ifNESW@\let\next@\lower\else\let\next@\raise\fi
     \next@\dimen@\hbox{\arrowfont@\char\angcount@}}}%
   \rlap{\hskip-\dimen@ii
    \smash{\ifNESW@\let\next@\raise\else\let\next@\lower\fi
      \next@\dimen@\hbox{\arrowfont@\char\angcount@}}}\hfil}%
 \fi
 \rlap{\smash{\hskip\tocenter@\hskip\firstx@
  \ifnum\arrcount@=\m@ne
  \else
   \ifnum\arrcount@=\thr@@
   \else
    \ifnum\scount@=\m@ne
    \else
     \ifnum\scount@=\z@
     \else
      \setbox\zer@\hbox{\ifnum\angcount@=117 \arrow@v\else\arrowfont@\fi
       \char\Scount@}%
      \ifNESW@
       \ifnum\scount@=\tw@
        \dimen@=\shifted@ \advance\dimen@-\charht@
        \ifN@\hskip-\wd\zer@\fi
        \Nnext@
        \next@\dimen@\copy\zer@
        \ifN@\else\hskip-\wd\zer@\fi
       \else
        \Nnext@
        \ifN@\else\hskip-\wd\zer@\fi
        \next@\shifted@\copy\zer@
        \ifN@\hskip-\wd\zer@\fi
       \fi
       \ifnum\scount@=12
        \advance\shifted@\charht@ \advance\goal@-\charht@
        \ifN@ \hskip\wd\zer@ \else \hskip-\wd\zer@ \fi
       \fi
       \ifnum\scount@=13
        \getcos@{\thr@@\p@}%
        \ifN@ \hskip\dimen@ \else \hskip-\wd\zer@ \hskip-\dimen@ \fi
        \adjust@\shifted@ \advance\adjust@\dimen@ii
        \Nnext@
        \next@\adjust@\copy\zer@
        \ifN@ \hskip-\dimen@ \hskip-\wd\zer@ \else \hskip\dimen@ \fi
       \fi
      \else
       \ifN@\hskip-\wd\zer@\fi
       \ifnum\scount@=\tw@
        \ifN@ \hskip\wd\zer@ \else \hskip-\wd\zer@ \fi
        \dimen@=\shifted@ \advance\dimen@-\charht@
        \Nnext@
        \next@\dimen@\copy\zer@
        \ifN@\hskip-\wd\zer@\fi
       \else
        \Nnext@
        \next@\shifted@\copy\zer@
        \ifN@\else\hskip-\wd\zer@\fi
       \fi
       \ifnum\scount@=12
        \advance\shifted@\charht@ \advance\goal@-\charht@
        \ifN@ \hskip-\wd\zer@ \else \hskip\wd\zer@ \fi
       \fi
       \ifnum\scount@=13
        \getcos@{\thr@@\p@}%
        \ifN@ \hskip-\wd\zer@ \hskip-\dimen@ \else \hskip\dimen@ \fi
        \adjust@\shifted@ \advance\adjust@\dimen@ii
        \Nnext@
        \next@\adjust@\copy\zer@
        \ifN@ \hskip\dimen@ \else \hskip-\dimen@ \hskip-\wd\zer@ \fi
       \fi      
      \fi
  \fi\fi\fi\fi
  \ifnum\arrcount@=\m@ne
  \else
   \loop
    \ifdim\goal@>\charht@
    \ifE@\else\hskip-\charwd@\fi
    \Nnext@
    \next@\shifted@\copy\shaft@
    \ifE@\else\hskip-\charwd@\fi
    \advance\shifted@\charht@ \advance\goal@ -\charht@
    \repeat
   \ifdim\goal@>\z@
    \dimen@=\charht@ \advance\dimen@-\goal@
    \divide\dimen@\tan@i \multiply\dimen@\tan@ii
    \ifE@ \hskip-\dimen@ \else \hskip-\charwd@ \hskip\dimen@ \fi
    \adjust@=\shifted@ \advance\adjust@-\charht@ \advance\adjust@\goal@
    \Nnext@
    \next@\adjust@\copy\shaft@
    \ifE@ \else \hskip-\charwd@ \fi
   \else
    \adjust@=\shifted@ \advance\adjust@-\charht@
   \fi
  \fi
  \ifnum\arrcount@=\m@ne
  \else
   \ifnum\arrcount@=\thr@@
   \else
    \ifnum\tcount@=\m@ne
    \else
     \setbox\zer@
      \hbox{\ifnum\angcount@=117 \arrow@v\else\arrowfont@\fi\char\Tcount@}%
     \ifnum\tcount@=\thr@@
      \advance\adjust@\charht@
      \ifE@\else\ifN@\hskip-\charwd@\else\hskip-\wd\zer@\fi\fi
     \else
      \ifnum\tcount@=12
       \advance\adjust@\charht@
       \ifE@\else\ifN@\hskip-\charwd@\else\hskip-\wd\zer@\fi\fi
      \else
       \ifE@\hskip-\wd\zer@\fi
     \fi\fi
     \Nnext@
     \next@\adjust@\copy\zer@
     \ifnum\tcount@=13
      \hskip-\wd\zer@
      \getcos@{\thr@@\p@}%
      \ifE@\hskip-\dimen@ \else\hskip\dimen@ \fi
      \advance\adjust@-\dimen@ii
      \Nnext@
      \next@\adjust@\box\zer@
     \fi
  \fi\fi\fi}}%
 \iflabel@i
  \rlap{\hskip\tocenter@
  \dimen@\firstx@ \advance\dimen@\secondx@ \divide\dimen@\tw@
  \advance\dimen@\ldimen@i
  \dimen@ii\firsty@ \advance\dimen@ii\secondy@ \divide\dimen@ii\tw@
  \multiply\ldimen@i\tan@i \divide\ldimen@i\tan@ii
  \ifNESW@ \advance\dimen@ii\ldimen@i \else \advance\dimen@ii-\ldimen@i \fi
  \setbox\zer@\hbox{\ifNESW@\else\ifnum\arrcount@=\thr@@\hskip4\p@\else
   \hskip\tw@\p@\fi\fi
   $\m@th\tsize@@\label@i$\ifNESW@\ifnum\arrcount@=\thr@@\hskip4\p@\else
   \hskip\tw@\p@\fi\fi}%
  \ifnum\arrcount@=\m@ne
   \ifNESW@ \advance\dimen@.5\wd\zer@ \advance\dimen@\p@ \else
    \advance\dimen@-.5\wd\zer@ \advance\dimen@-\p@ \fi
   \advance\dimen@ii-.5\ht\zer@
  \else
   \advance\dimen@ii\dp\zer@
   \ifnum\slcount@<6 \advance\dimen@ii\tw@\p@ \fi
  \fi
  \hskip\dimen@
  \ifNESW@ \let\next@\llap \else\let\next@\rlap \fi
  \next@{\smash{\raise\dimen@ii\box\zer@}}}%
 \fi
 \iflabel@ii
  \ifnum\arrcount@=\m@ne
  \else
   \rlap{\hskip\tocenter@
   \dimen@\firstx@ \advance\dimen@\secondx@ \divide\dimen@\tw@
   \ifNESW@ \advance\dimen@\ldimen@ii \else \advance\dimen@-\ldimen@ii \fi
   \dimen@ii\firsty@ \advance\dimen@ii\secondy@ \divide\dimen@ii\tw@
   \multiply\ldimen@ii\tan@i \divide\ldimen@ii\tan@ii
   \advance\dimen@ii\ldimen@ii
   \setbox\zer@\hbox{\ifNESW@\ifnum\arrcount@=\thr@@\hskip4\p@\else
    \hskip\tw@\p@\fi\fi
    $\m@th\tsize@\label@ii$\ifNESW@\else\ifnum\arrcount@=\thr@@\hskip4\p@
    \else\hskip\tw@\p@\fi\fi}%
   \advance\dimen@ii-\ht\zer@
   \ifnum\slcount@<9 \advance\dimen@ii-\thr@@\p@ \fi
   \ifNESW@ \let\next@\rlap \else \let\next@\llap \fi
   \hskip\dimen@\next@{\smash{\raise\dimen@ii\box\zer@}}}%
  \fi
 \fi
}
\def\outnewCD@#1{\def#1{\Err@{\string#1 must not be used within \string\newCD}}}
\newskip\prenewCDskip@
\newskip\postnewCDskip@
\prenewCDskip@\z@
\postnewCDskip@\z@
\def\prenewCDspace#1{\RIfMIfI@
 \onlydmatherr@\prenewCDspace\else\advance\prenewCDskip@#1\relax\fi\else
 \onlydmatherr@\prenewCDspace\fi}
\def\postnewCDspace#1{\RIfMIfI@
 \onlydmatherr@\postnewCDspace\else\advance\postnewCDskip@#1\relax\fi\else
 \onlydmatherr@\postnewCDspace\fi}
\def\predisplayspace#1{\RIfMIfI@
 \onlydmatherr@\predisplayspace\else
 \advance\abovedisplayskip#1\relax
 \advance\abovedisplayshortskip#1\relax\fi
 \else\onlydmatherr@\prenewCDspace\fi}
\def\postdisplayspace#1{\RIfMIfI@
 \onlydmatherr@\postdisplayspace\else
 \advance\belowdisplayskip#1\relax
 \advance\belowdisplayshortskip#1\relax\fi
 \else\onlydmatherr@\postdisplayspace\fi}
\def\PrenewCDSpace#1{\global\prenewCDskip@#1\relax}
\def\PostnewCDSpace#1{\global\postnewCDskip@#1\relax}
\def\newCD#1\endnewCD{%
 \outnewCD@\cgaps\outnewCD@\rgaps\outnewCD@\Cgaps\outnewCD@\Rgaps
 \prenewCD@#1\endnewCD
 \advance\abovedisplayskip\prenewCDskip@
 \advance\abovedisplayshortskip\prenewCDskip@
 \advance\belowdisplayskip\postnewCDskip@
 \advance\belowdisplayshortskip\postnewCDskip@
 \vcenter{\vskip\prenewCDskip@ \Let@ \colcount@\@ne \rowcount@\z@
  \everycr{%
   \noalign{%
    \ifnum\rowcount@=\Rowcount@
    \else
     \global\nointerlineskip
     \getrgap@\rowcount@ \vskip\getdim@
     \global\advance\rowcount@\@ne \global\colcount@\@ne
    \fi}}%
  \tabskip\z@
  \halign{&\global\xoff@\z@ \global\yoff@\z@
   \getcgap@\colcount@ \hskip\getdim@
   \hfil\vrule height10\p@ width\z@ depth\z@
   $\m@th\displaystyle{##}$\hfil
   \global\advance\colcount@\@ne\cr
   #1\crcr}\vskip\postnewCDskip@}%
 \prenewCDskip@\z@\postnewCDskip@\z@
 \def\getcgap@##1{\ifcase##1\or\getdim@\z@\else\getdim@\standardcgap\fi}%
 \def\getrgap@##1{\ifcase##1\getdim@\z@\else\getdim@\standardrgap\fi}%
 \let\Width@\relax\let\Height@\relax\let\Depth@\relax\let\Rowheight@\relax
 \let\Rowdepth@\relax\let\Colwdith@\relax
}
\catcode`\@=\active

\hsize 30pc
\vsize 47pc
\magnification=\magstep1
\let\[\lceil
\let\]\rfloor
\def\nmb#1#2{#2}         
\def\cit#1#2{\ifx#1!\cite{#2}\else#2\fi} 
\def\idx{}               
\def\ign#1{}             
\redefine\o{\circ}
\let\ceylon\colon
\redefine\colon{\hskip.05em\ceylon}
\define\({\big(}
\define\){\big)}
\define\X{\frak X}
\define\al{\alpha}
\define\be{\beta}

\define\de{\delta}

\define\la{\lambda}
\define\rh{\rho}
\define\si{\sigma}
\define\ta{\tau}
\define\ph{\varphi}

\define\ps{\psi}

\define\Ph{\Phi}
\define\Ps{\Psi}

\define\x{\times}
\define\Id{\operatorname{Id}}
\define\g{\frak g}
\define\h{\frak h}
\define\e{\frak e}

\redefine\l{\frak l}
\define\ad{\operatorname{ad}}
\define\pr{\operatorname{pr}}

\define\der{\operatorname{der}}

\define\s#1#2{\operatorname{sign}(#1,\bold #2)}
\redefine\L{\operatorname{\Cal L}}
\long\def\alert#1{\par\medskip{\narrower%
     \vrule\vbox{\advance\hsize by -2.1\parindent\hrule
     \medskip\noindent#1\medskip\hrule}\vrule}\medskip}

\def\today{\ifcase\month\or
 January\or February\or March\or April\or May\or June\or
 July\or August\or September\or October\or November\or December\fi
 \space\number\day, \number\year}
\topmatter
\title  Extensions of super Lie algebras 
\endtitle
\author Dmitri Alekseevsky, Peter W. Michor, Wolfgang Ruppert  
\endauthor
\affil
Erwin Schr\"odinger Institut f\"ur Mathematische Physik,\endgraf
Boltzmanngasse 9, A-1090 Wien, Austria
\endaffil
\address 
D.V. Alekseevsky: 
Department of Mathematics,
University of Hull,
Cottingham Road,
Hull, HU6 7RX,
England
\endaddress
\email
d.v.alekseevsky\@maths.hull.ac.uk 
\endemail
\address
P\. W\. Michor: Institut f\"ur Mathematik, Universit\"at Wien,
Strudlhofgasse 4, A-1090 Wien, Austria; {\it and:} 
Erwin Schr\"odinger Institut f\"ur Mathematische Physik,
Boltzmanngasse 9, A-1090 Wien, Austria
\endaddress
\email Peter.Michor\@esi.ac.at \endemail
\address
W\. Ruppert: 
Institut f\"ur Mathematik und angewandte Statistik, Universit\"at 
f\"ur Bodenkultur,
Gregor Mendelstrasse 33, A-1180 Wien, Austria
\endaddress
\email Ruppert\@edv1.boku.ac.at \endemail

\dedicatory \enddedicatory
\date {\today} \enddate
\thanks  
P.W.M. was supported  
by `Fonds zur F\"orderung der wissenschaftlichen  
Forschung, Projekt P~14195~MAT'.
\endthanks 
\keywords Extensions of super Lie algebras, cohomology of super Lie algebras 
\endkeywords
\subjclass\nofrills{\rm 2000}
 {\it Mathematics Subject Classification}.\usualspace
 Primary 17B05, 17B56\endsubjclass

\abstract We study (non-abelian) extensions of a given super Lie algebra, 
identify a cohomological obstruction to the existence, parallel to the
known one for Lie algebras.
An analogy to the setting of covariant exterior derivatives, curvature, 
and the Bianchi identity in differential geometry is spelled out. 
\endabstract
\endtopmatter

\long\def\alert#1{\relax}
\document

\subhead\nmb0{1}. Introduction \endsubhead
The theory of group extensions and their interpretation in terms of 
cohomology is  well known, see, e.g., \cit!{3}, \cit!{6}, 
\cit!{4}, \cit!{2}.
Analogous results for Lie algebras are dispersed in the literature, see
\cit!{5}, \cit!{15}, \cit!{19}; 
the case of Lie algebroids is treated in \cit!{14}. 
We owe this to Kirill Mackenzie.

This paper gives a unified and coherent account of this for 
super Lie algebras stressing a certain analogy 
with concepts from differntial geometry:
Covariant exterior 
derivatives, curvature and the Bianchi identity.

In an unpublished preliminary version of this paper \cit!{1}, 
the analogous results for Lie algebras were developed. 

\subhead\nmb0{2}. Super Lie algebras \endsubhead
(See \cit!{8}, or \cit!{16} for an introduction) 
A {\it super Lie algebra} is a 2-graded vector space 
$\g=\g_0\oplus \g_1$, together with a graded  Lie bracket 
$[\quad,\quad]:\g\x \g\to \g$ of degree 0. That is, $[\ ,\ ]$ is a bilinear 
map with 
$[\g_i,\g_j]\subseteq \g_{i+j (\text{mod} 2)}$, and such that for 
homogeneous elements $X\in \g_x$, $Y\in \g_y$, and $Z\in \g_z$ the 
identities   
$$\alignat2
[X,Y]&=-(-1)^{xy}[Y,X] &\quad&\text{ (graded antisymmetry) }\\
[X,[Y,Z]] &= [[X,Y],Z] +(-1)^{xy}[Y,[X,Z]] &&\text{ (graded Jacobi identity)}
\endalignat$$
hold.
The graded Jacobi identity, shorter 
$\sum_{\text{cyclic}} (-1)^{xz} [X,[Y,Z]] =0$,
says that $\ad_X:\g\to\g, \ Y\mapsto [X,Y]$ is a graded 
derivation of degree $x$, so that 
$\ad_X[Y,Z]=[\ad_XY,Z]+(-1)^{xy}[X,\ad_XZ]$.
We denote by $\der(\g)$ the super Lie algebra of graded derivations 
of $\g$. The notion of homomorphism is as usual, homomorphisms are 
always of degree 0.  

\subhead\nmb0{3}. Describing extensions, first part \endsubhead
Consider any exact sequence of homomorphisms of super Lie algebras: 
$$
0\to \h @>i>> \e @>p>> \g \to 0
$$
Consider a graded linear mapping $s:\g \to \e$ of degree 0 with 
$p\o s=\Id_\g$. Then $s$ induces mappings 
$$\gather
\al:\g\to \der(\h)
     \quad \text{({\it super connection}) by }\quad 
     \al_X(H)= [s(X),H], \tag{\nmb.{3.1}}\\
\rh:\bigwedge^2_{\text{graded}}\g \to \h
     \quad \text{({\it curvature}) by }\quad 
     \rh(X,Y)= [s(X),s(Y)]-s([X,Y]) \tag{\nmb.{3.2}}\\
\endgather$$
which are easily seen to be of degree 0 and to satisfy: 
$$\gather
[\al_X,\al_Y]-\al_{[X,Y]} = \ad_{\rh(X,Y)} \tag{\nmb.{3.3}}\\
\sum_{\text{cyclic}} (-1)^{xz} \Bigl(\al_X \rh(Y,Z) 
     - \rh([X,Y],Z) \Bigr) =0\tag{\nmb.{3.4}}
\endgather$$
Property \thetag{\nmb!{3.4}} is equivalent to the graded Jacobi identity in 
$\e$. 

\subhead\nmb0{4}. Motivation: Lie algebra extensions associated to a 
principal bundle
\endsubhead
In the case of Lie algebras, the extension $$
0\to \h @>i>> \e @>p>> \g \to 0
$$
appears in the following geometric situation. Let $\pi : P \to M = 
P/K $ be a principal bundle with structure group $K$. Then the Lie 
algebra of infinitesimal automorphisms $\e=\X(P)^K$, i.e.\ the Lie 
algebra of $K$-invariant vector fields on $P$, is an extension of the 
Lie algebra $\g=\X(M)$ of all vector fields on $M$ by the Lie algebra 
$\h=\X_{\text{vert}}(P)^K$ of all vertical $K$-invariant vector 
fields, i.e., infinitesimal gauge transformations. 
In this case we have simultaneously an extension of 
$C^\infty(M)$-modules. A section $s:\g\to\e$ which is simultaneously 
a homomorphism of $C^\infty(M)$-modules can be considered as a 
connection, and $\rh$, defined as in \nmb!{3.2}, is the curvature of 
this connection. This geometric example is a guideline for our 
approach. It works also for super Lie algebras. See \cit!{9}, 
section~11 for more background information. This analogy with 
differential geometry has also been noticed in \cit!{10} and 
\cit!{11} and has been used used extensively in the theory of Lie 
algebroids, see \cit!{14}. 

\subhead\nmb0{5}. Algebraic theory of connections, curvature, and 
cohomology \endsubhead
We want to interpret \nmb!{3.4} as $\de_\al\rh=0$ 
where $\de_\al$ is 
an analogon of the graded version of the Chevalley coboundary 
operator, but with values in the non-representation $\h$; we shall 
see that this is exactly the notion of a 
\idx{\it super exterior covariant derivative}. Namely, let 
$L^{p,y}_{\text{gskew}}(\g;\h)$ be the space of all graded 
antisymmetric $p$-linear 
mappings $\Ph:\g^p\to \h$ of degree $y$, i.e\. 
$$\gather
\Ph(X_1,\dots,X_p)\in \h_{y+x_1+\dots+x_p}, \\
\Ph(X_1,\dots,X_p)=-(-1)^{x_ix_{i+1}}\Ph(X_1,\dots,X_{i+1},X_i,\dots,X_p).
\endgather$$ 
In order to treat the graded Chevalley coboundary operator we 
need the following notation, which is similar to the one used in 
\cit!{12},~3.1:
Let $\bold x = (x_1,\dots,x_k) \in (\Bbb Z_2)^k$ be a multi index of 
binary degrees
$x_i \in \Bbb Z_2$ and let $\si \in \Cal S_k$ be a
permutation of $k$ symbols. Then we define the {\it multigraded
sign} $\s \si x$ as follows: 
For a transposition $\si =(i,i+1)$ we
put $\s \si x = -(-1)^{x_i\,x_{i+1}}$; 
it can be checked by combinatorics that this gives a well defined
mapping $\s {\quad}x:\Cal S_k \to  \{-1,+1\}$. In fact one may define
directly
$\s \si x = \operatorname{sign}(\si)
\operatorname{sign}(\si_{|x_1|,\dots,|x_k|})$,
where $|\quad|: \Bbb Z_2\to \Bbb Z$ is the embedding and  
where $\si_{|x_1|,\ldots,|x_k|}$ is that permutation of 
$|x_1|+\dots+|x_k|$ symbols which moves the $i$-th block of
length $|x_i|$ to the position $\si i$, and where 
$\operatorname{sign}(\si)$ 
denotes the ordinary sign of a permutation in $\Cal S_k$.
Let us write $\si x = (x_{\si1},\dots,x_{\si k})$, 
then we have 
$$
\s {\si\o\ta}x = \s \si x .\s \ta{\si \bold x},
$$
and $\Ph\in L^{p,y}_{\text{gskew}}(\g;\h)$ satisfies
$$
\Ph(X_{\si1},\dots,X_{\si p}) = \s \si x \Ph(X_1,\dots,X_p)
$$
Given a super connection $\al:\g\to \der(\h)$ as in \nmb!{3.1},
we define the graded version of the covariant exterior derivative
by  
$$\align
\de_\al&:L^{p,y}_{\text{gskew}}(\g;\h) \to L^{p+1,y}_{\text{gskew}}(\g;\h)\\
(\de_\al\Ph)(X_0,\dots,X_p) 
&= \sum_{i=0}^p(-1)^{x_iy+a_i(\bold x)}\al_{X_i}
     (\Ph(X_0,\dots,\widehat{X_i},\dots,X_p))\\
&\quad+ \sum_{i<j}(-1)^{a_{ij}(\bold x)}
     \Ph([X_i,X_j],X_0,\dots,\widehat{X_i},\dots,
     \widehat{X_j},\ldots,X_p)\\
a_i(\bold x) &= x_i(x_1+\dots+x_{i-1}) + i\\
a_{ij}(\bold x) &= a_i(\bold x) + a_j(\bold x) + x_ix_j
\endalign$$
for cochains $\Ph$ with coefficients in the non-representation $\h$ 
of $\g$. In fact,  
$\de_\al$ has the formal property of a 
\idx{\it super covariant exterior derivative}, namely:
$$
\de_\al(\ps\wedge \Ph) = \de\ps \wedge \Ph + 
(-1)^{q}\ps\wedge \de_\al\Ph
$$
for $\Ph\in L^{p,y}_{\text{gskew}}(\g;\h)$ and 
$\ps\in L^{q,z}_{\text{gskew}}(\g;\Bbb R)$ a form of degree $q$ and 
weight $z$ (we put $\Bbb R$ of degree 0), where 
$$
(\de\ps)(X_0,\dots,X_q) 
= \sum_{i<j}(-1)^{a_{ij}(\bold x)}
     \Ph([X_i,X_j],X_0,\dots,\widehat{X_i},\dots,
     \widehat{X_j},\ldots,X_q)
$$
is the super analogon of the Chevalley coboundary operator for 
cochains with values in the trivial $\g$-representation $\Bbb R$, and 
where the module structure is given by
$$\multline
(\ps\wedge \Ph)(X_1\dots,X_{q+p}) = \\
= \frac1{q!p!}\sum_{\si\in \Cal S_{q+p}} \s \si x 
(-1)^{yb_q(\si,\bold x)} 
\ps(X_{\si1},\dots,X_{\si q})\Ph(X_{\si(q+1)},\dots,X_{\si(q+p)}),
\endmultline$$
where $b_i(\si,\bold x)=|x_{\si1}|+\dots+|x_{\si i}|$.

Moreover for $\Ph\in L^{p,y}_{\text{gskew}}(\g;\h)$ and 
$\Ps\in L^{q,z}_{\text{gskew}}(\g;\h)$ we put
$$\multline
[\Ph,\Ps]_\wedge(X_1,\dots,X_{p+q}) = \\
= \frac 1{p!\,q!} \sum_{\si} \s\si x (-1)^{zb_p(\si,\bold x)}
        [\Ph(X_{\si1},\dots,X_{\si p}),
     \Ps(X_{\si(p+1)},\dots,X_{\si(p+q)})]_{\h}.
\endmultline$$

\noindent \therosteritem{\nmb.{5.1}} The bracket
$[\quad,\quad]_\wedge$ is a $\Bbb Z\x \Bbb Z_2$-graded Lie algebra 
structure on  
$$
L^*_{\text{skew}}(V,\h)=\bigoplus_{p\in \Bbb Z_{\ge0},y\in \Bbb Z_2} 
L^{p,y}_{\text{gskew}}(\g;\h)
$$
which means that the analoga of the properties of section \nmb!{2} 
hold for the signs $(-1)^{p_1p_2+y_1y_2}$. See \cit!{12} for more 
details.  

A straightforward computation shows that
for $\Ph\in L^{p,y}_{\text{gskew}}(\g;\h)$ we have
$$\de_\al\de_\al(\Ph) = [\rh,\Ph]_\wedge.\tag{\nmb.{5.2}}$$

Note that \nmb!{5.2} justifies the use of the super analogon of the 
Chevalley cohomology if $\al:\g\to\der(\h)$ 
is a homomorphism of super 
Lie algebras or $\al:\g\to \operatorname{End}(V)$ is a representation 
in a graded vector space. 
See \cit!{12} for more details. 

\subhead\nmb0{6}. Describing extensions, continued \endsubhead
Continuing the discussion of section \nmb!{3},
we now can describe completely the super Lie algebra structure on 
$\e=\h\oplus s(\g)$ in terms of $\al$ and $\rh$:
$$\multline
[H_1+s(X_1),H_2+s(X_2)] = \\ 
=([H_1,H_2]+\al_{X_1}H_2 -(-1)^{h_1x_2}\al_{X_2}H_1 
+\rh(X_1,X_2)) + s[X_1,X_2]. 
\endmultline\tag{\nmb.{6.1}}$$
If $\al:\g\to\der(\h)$ and 
$\rh:\bigwedge^2_{\text{graded}}\g\to \h$ satisfy \thetag{\nmb!{3.3}} and 
\thetag{\nmb!{3.4}} then one checks easily that formula 
\thetag{{\nmb!{6.1}}} gives a super Lie algebra structure  
on $\h\oplus s(\g)$. 

If we change the linear section $s$ to $s'=s+b$ for linear 
$b:\g\to\h$ of degree zero, then we get 
$$\align
\al'_X &= \al_X + \ad^\h_{b(X)} \tag{\nmb.{6.2}}\\
\rh'(X,Y) &= \rh(X,Y) + \al_Xb(Y) -(-1)^{xy}\al_Yb(X) - b([X,Y]) + [bX,bY] 
\tag{\nmb.{6.3}}\\
&= \rh(X,Y) + (\de_\al b)(X,Y) + [bX,bY].\\
\rh'&= \rh + \de_\al b + \tfrac12 [b,b]_\wedge. 
\endalign$$

\proclaim{\nmb0{7}. Proposition}
Let $\h$ and $\g$ be super Lie algebras. 

Then isomorphism classes of extensions of $\g$ over $\h$, i.e\. short 
exact sequences of Lie algebras 
$0\to \h\to \e\to \g\to 0$ modulo the equivalence described by 
commutative diagrams of super Lie algebra homomorphisms
$$\CD
0 @>>> \h @>>> \e @>>> \g @>>> 0 \\
@.     @|  @V{\ph}VV   @|  @. \\
0 @>>> \h @>>> \e' @>>> \g @>>> 0 \\
\endCD$$
correspond bijectively
to equivalence classes of data of the following form:
\smallskip\noindent
\thetag{\nmb.{7.1}} \quad a linear mapping 
     $\al:\g\to \der(\h)$ of degree 0,
\smallskip\noindent
\thetag{\nmb.{7.2}} \quad  a graded skew-symmetric bilinear mapping 
     $\rh:\g\x\g\to \h$ of degree 0, 
\smallskip\noindent
such that
$$\align
&[\al_X,\al_Y]-\al_{[X,Y]} = \ad_{\rh(X,Y)} \tag{\nmb.{7.3}}\\
&\sum_{\text{cyclic}} (-1)^{xz}\Bigl(\al_X \rh(Y,Z) - \rh([X,Y],Z) \Bigr) 
=0.\tag{\nmb.{7.4}}
\endalign$$
On the vector space $\e:=\h\oplus \g$ a Lie algebra structure is 
given by 
$$\multline
[H_1+X_1,H_2+X_2]_\e = \\
= ([H_1,H_2]_\h + \al_{X_1}H_2 
     -(-1)^{x_2h_1} \al_{X_2}H_1 + \rh(X_1,X_2)) + 
     [X_1,X_2]_\g,
\endmultline\tag{\nmb.{7.5}}$$
and the associated exact sequence is  
$$
0\to \h @>i_1>>\h\oplus\g= \e @>{\operatorname{pr}_2}>> \g\to 0.$$
Two data $(\al,\rh)$ and $(\al',\rh')$ are equivalent if there exists
a linear mapping $b:\g\to\h$ of degree 0 such that 
$$\align
\al'_X &= \al_X + \ad^\h_{b(X)}, \tag{\nmb.{7.6}}\\
\rh'(X,Y) &= \rh(X,Y) + \al_Xb(Y) 
     -(-1)^{xy}\al_Yb(X) - b([X,Y]) + [b(X),b(Y)] \\
\rh'&= \rh + \de_\al b + \tfrac12 [b,b]_\wedge, \tag{\nmb.{7.7}}\\
\endalign$$
the corresponding isomorphism being
$$
\e=\h\oplus\g \to \h\oplus\g=\e',\qquad H+X\mapsto H-b(X)+ X.
$$
Moreover, a datum $(\al,\rh)$ corresponds to a split extension (a 
semidirect product) if and only if $(\al,\rh)$ is equivalent to to a 
datum of the form $(\al',0)$ (then $\al'$ is a homomorphism). This is 
the case if and only if there exists a mapping $b:\g\to \h$ such that
$$
\rh = \de_\al b - \tfrac12[b,b]_\wedge.\tag{\nmb.{7.8}} 
$$
\endproclaim
\demo{Proof}
Direct computations. 
\qed\enddemo

\proclaim{\nmb0{8}. Corollary}
Let $\g$ and $\h$ be super Lie algebras such that $\h$ has no 
(graded) center. 
Then isomorphism classes of extensions of $\g$ over $\h$ correspond 
bijectively to homomorphisms of super Lie algebras
$$
\bar\al:\g\to \operatorname{out}(\h)=\der(\h)/\ad(\h). 
$$
\endproclaim

\demo{Proof} Choose a linear lift $\al:\g\to \der(\h)$ of 
$\bar\al$.
Since $\bar\al:\g\to \der(\h)/\ad(\h)$ is a homomorphism, 
there is a uniquely defined skew symmetric linear mapping
$\rh:\g\x\g\to \h$ such that
$[\al_X,\al_Y]-\al_{[X,Y]} = \ad_{\rh(X,Y)}$.
Condition \thetag{\nmb!{7.4}} is then automatically 
satified. For later use we record the simple proof:
$$\align
&\sum_{\text{cyclic} X,Y,Z}(-1)^{xz}\Bigl[\al_X \rh(Y,Z) 
     - \rh([X,Y],Z),H \Bigr] \\
&=\sum_{\text{cyclic} X,Y,Z}(-1)^{xz}\Bigl(\al_X [\rh(Y,Z),H] 
     -(-1)^{(x(y+z))} [\rh(Y,Z),\al_XH] -\\
&\qquad\qquad\qquad\qquad\qquad\qquad\qquad\qquad\qquad\qquad
- [\rh([X,Y],Z),H] \Bigr) \\
&=\sum_{\text{cyclic} X,Y,Z}(-1)^{xz}\Bigl(\al_X [\al_Y,\al_Z] 
     - \al_X\al_{[Y,Z]} -(-1)^{(x(y+z))} [\al_Y,\al_Z]\al_X +\\
&\qquad\qquad\qquad\qquad\qquad
     +(-1)^{(x(y+z))} \al_{[Y,Z]}\al_X 
     - [\al_{[X,Y]},\al_Z] + \al_{[[X,Y]Z]} \Bigr)H \\
&=\sum_{\text{cyclic} X,Y,Z}(-1)^{xz}\Bigl([\al_X,[\al_Y,\al_Z]] 
     -[\al_X,\al_{[Y,Z]}] - [\al_{[X,Y]},\al_Z] 
     + \al_{[[X,Y]Z]} \Bigr)H =0.\\
\endalign$$
Thus $(\al,\rh)$ describes an extension, by \nmb!{7}. 
The rest is clear.
\qed\enddemo

\subhead\nmb0{9}. Remark \endsubhead
If the super Lie algebra $\h$ has no center and a homomorphism 
$\bar\al:\g\to \operatorname{out}(\h)=\der(\h)/\ad(\h)$ is given, the 
extension corresponding to $\bar\al$ 
is given by the pullback diagram
$$\CD
0 @>>>\h @>>>\der(\h) \x_{\operatorname{out}(\h)} \g @>{\pr_2}>>\g @>>> 0\\
@.     @|       @V{\pr_1}VV                            @V{\bar\al}VV @.\\
0 @>>> \h @>>> \der(\h) @>{\pi}>> \operatorname{out}(\h) @>>> 0\\
\endCD$$
where $\der(\h) \x_{\operatorname{out}(\h)} \g$ is the Lie subalgebra
$$
\der(\h) \x_{\operatorname{out}(\h)} \g 
     := \{(D,X)\in \der(\h)\x \g:\pi(D)=\bar\al(X)\} \subset 
     \der(\h)\x\g.
$$
We owe this remark to E\. Vinberg.

If the super Lie algebra $\h$ has no center and satisfies 
$\der(\h)=\h$, and if $\h$ is  
an ideal in a super Lie algebra $\e$, then $\e\cong \h\oplus\e/\h$, since 
$\operatorname{Out}(\h)=0$.

\proclaim{\nmb0{10}. Theorem}
Let $\g$ and $\h$ be super Lie algebras and let 
$$\bar\al:\g\to \operatorname{out}(\h)=\der(\h)/\ad(\h)$$
be a homomorphism of super Lie algebras. 
Then the following are equivalent:
\roster
\item"(\nmb.{10.1})" For one (equivalently: any) linear lift 
       $\al:\g\to \der(\h)$ of degree 0 of $\bar\al$ choose 
       $\rh:\bigwedge^2_{\text{graded}}\g\to \h$ of degree 0 satisfying 
       $([\al_X,\al_Y]-\al_{[X,Y]})=\ad_{\rh(X,Y)}$. Then the 
       $\de_{\bar\al}$-cohomology class of 
       $\la=\la(\al,\rh):=\de_{\al}\rh:\bigwedge^3\g \to Z(\h)$ in 
       $H^3(\g;Z(\h))$ vanishes. 
\item "(\nmb.{10.2})" There exists an extension 
       $0\to \h\to\e\to\g\to 0$ inducing the homomorphism $\bar\al$. 
\endroster
If this is the case then all extensions $0\to \h\to\e\to\g\to 0$ 
inducing the homomorphism $\bar\al$ are parameterized by 
$H^2(\g,(Z(\h),\bar\al))$, the second graded Chevalley cohomology space of 
the super Lie algebra $\g$ with values in the graded $\g$-module 
$(Z(\h),\bar\al)$. 
\endproclaim
\demo{Proof}
From the computation in the proof of corollary \nmb!{8} it follows that 
$$
\ad(\la(X,Y,Z))=\ad(\de_{\al}\rh(X,Y,Z))=0
$$
so that 
$\la(X,Y,Z)\in Z(\h)$.
The super Lie algebra $\operatorname{out}(\h)=\der(\h)/\ad(\h)$ 
acts on the center $Z(\h)$, thus $Z(\h)$ is a graded $\g$-module via 
$\bar \al$, and $\de_{\bar\al}$ is the differential of the Chevalley 
cohomology. 
Using \nmb!{5.2}, then \nmb!{5.1} we see
$$
\de_{\bar\al}\la = \de_{\al}\de_{\al}\rh = 
     [\rh,\rh]_\wedge = -(-1)^{2\cdot 2+ 0\cdot 0}[\rh,\rh]_\wedge =0,
$$
so that $[\la]\in H^3(\g;Z(\h))$. 

Let us check next that the cohomology class $[\la]$ does not depend 
on the choices we made. If we are given a pair $(\al,\rh)$ as above 
and we take another linear lift $\al':\g\to\der(\h)$ then 
$\al'_X=\al_X+\ad_{b(X)}$ for some linear $b:\g\to\h$. We consider 
$$
\rh':\bigwedge^2_{\text{graded}}\g\to\h,\quad 
     \rh'(X,Y)=\rh(X,Y)+(\de_\al b)(X,Y)+[b(X),b(Y)]. 
$$
Easy computations show that
$$\gather
[\al'_X,\al'_Y]-\al'_{[X,Y]} = \ad_{\rh'(X,Y)}\\
\la(\al,\rh) = \de_\al\rh = \de_{\al'}\rh' = \la(\al',\rh')
\endgather$$ 
so that even the cochain did not change. 
So let us consider for fixed $\al$ two linear mappings 
$$
\rh,\rh':\bigwedge^2_{\text{graded}}\g\to \h,\quad 
[\al_X,\al_Y]-\al_{[X,Y]} = \ad_{\rh(X,Y)}= \ad_{\rh'(X,Y)}. 
$$
Then $\rh-\rh'=:\mu:\bigwedge^2_{\text{graded}}\g\to Z(\h)$ and clearly
$\la(\al,\rh)-\la(\al,\rh')=\de_\al\rh-\de_\al\rh'=\de_{\bar\al}\mu$.

If there exists an extension inducing $\bar\al$ then for any lift 
$\al$ we may find $\rh$ as in proposition~\nmb!{7} such that 
$\la(\al,\rh)=0$. 
On the other hand, given a pair $(\al,\rh)$ as in \therosteritem1 
such that $[\la(\al,\rh)]=0\in H^3(\g,(Z(\h),\bar\al))$, there 
exists $\mu:\bigwedge^2\g\to Z(\h)$ such that $\de_{\bar\al}\mu=\la$.
But then 
$$
\ad_{(\rh-\mu)(X,Y)}=\ad_{\rh(X,Y)},\quad \de_\al(\rh-\mu)=0,
$$ 
so that $(\al,\rh-\mu)$ satisfy the conditions of \nmb!{7} and thus 
define an extension which induces $\bar\al$.

Finally, suppose that \therosteritem{\nmb!{10.1}} is satisfied, and 
let us determine how many extensions there exist which induce 
$\bar\al$. 
By proposition \nmb!{7} we have to determine all equivalence classes 
of data $(\al,\rh)$ as described there.
We may fix the linear lift $\al$ and one mapping 
$\rh:\bigwedge^2_{\text{graded}}\g\to \h$ which satisfies 
\thetag{\nmb!{7.3}} and \thetag{\nmb!{7.4}}, and we have to find all 
$\rh'$ with this property. But then 
$\rh-\rh'=\mu:\bigwedge^2_{\text{graded}}\g\to Z(\h)$ and 
$$
\de_{\bar\al}\mu = \de_\al\rh-\de_\al\rh'=0-0=0
$$
so that $\mu$ is a 2-cocycle. Moreover we may still 
pass to equivalent data in the sense of proposition~\nmb!{7} using some 
$b:\g\to\h$ which does not change $\al$, i.e\. $b:\g\to Z(\h)$. The 
corresponding $\rh'$ is, by \thetag{\nmb!{7.7}},
$\rh'=\rh+\de_\al b + \tfrac12[b,b]_{\wedge }=\rh+\de_{\bar\al}b$.
Thus only the cohomology class of $\mu$ matters.  
\qed\enddemo

\proclaim{\nmb0{11}. Corollary}
Let $\g$ and $\h$ be super Lie algebras such that $\h$ is abelian.
Then isomorphism classes of extensions of $\g$ over $\h$ correspond 
bijectively to the set of all pairs $(\al,[\rh])$, where
$\al:\g\to \g\frak l(\h)=\der(\h)$ is a homomorphism of super Lie algebras 
and $[\rh]\in H^2(\g,\h)$ is a graded Chevalley cohomology class with 
coefficients in the $\g$-module $\h$. 
\endproclaim

\demo{Proof} This is obvious from theorem \nmb!{10}. 
\qed\enddemo

\subhead\nmb0{12}. An interpretation of the class $\la$ \endsubhead
Let $\h$ and $\g$ be super Lie algebras and let a homomorphism of 
super Lie algebras 
$\bar\al:\g\to \der(\h)/\ad(\h)$ be given. We consider the extension 
$$
0\to \ad(\h) \to \der(\h) \to \der(\h)/\ad(\h) \to 0
$$
and the following diagram, where the bottom right hand square is a 
pullback (compare with remark \nmb!{9}):
$$\cgaps{0.5;0.5;0.5;1;1;0.5}\rgaps{0.5;1;1;0.4;1.2}\newCD
&&         & 0 @(0,-1) & 0 @(0,-1) & & \\
&&         & Z(\h) @()\a=@(1,0) @(0,-1) & Z(\h) @()\a-@(0,-1)  &    &  \\
&& 0 @(1,0) & \h @()\a-@(1,0) @(0,-1)    
     & \e @()\a-@(0,-1) @()\a-@(1,0) &  \g @(1,0) @()\a=@(0,-1) & 0  \\
&& 0 @(1,0) 
   & \ad(\h) @()\L{i}@(1,0) @()\a=\ds(4;0)\dtX(-4;0)@(-2,-2) @(0,-1) 
     & \tilde\e @()\L{\be}\l{\text{\quad pull back}}@(-1,-2) @()\L{p}@(1,0) @(0,-1)
       & \g @(1,0) @()\dtX(-10;0)\l{\bar\al}@(-1,-2) & 0\\
&&         & 0 & 0 &  & \\
0 @(1,0) & \ad(\h) @(2,0) &  
     & \der(\h) @(1,0) 
             &  \der(\h)/\ad(\h) @(1,0) & 0 & \\
\endnewCD$$

The left hand vertical column describes $\h$ as a central extension of 
$\ad(\h)$ with abelian kernel $Z(\h)$ which is moreover killed 
under the action of $\g$ via $\bar\al$; it is given by a cohomology 
class $[\nu]\in H^2(\ad(\h);Z(\h))^\g$. 
In order to get an extension $\e$ of $\g$  
with kernel $\h$ as in the third row we have to check that the 
cohomology class $[\nu]$ is in the image of 
$i^*:H^2(\tilde\e;Z(\h))\to H^2(\ad(\h);Z(\h))^\g$. 
It would be interesting to interpret this in terms of the super analogon of
the Hochschild-Serre spectral sequence from 
\cit!{7}.


\Refs

\widestnumber\key{ABCD}

\ref
\key \cit0{1} 
\by Alekseevsky, D.; Michor, Peter W.; Ruppert, Wolfgang A.F.
\paper Extensions of Lie algebras
\paperinfo unpublished
\finalinfo arXiv:math.DG/0005042
\endref

\ref
\key \cit0{2} 
\by Azc\'arraga, Jos\'e A.; Izquierdo, Jos\'e M.
\book Lie groups, Lie algebras, cohomology and some applications in 
physics
\bookinfo Cambridge Monographs on Mathematical Physics
\publ Cambridge University Press
\publaddr Cambridge, UK
\yr 1995
\endref

\ref
\key \cit0{3}
\by Eilenberg, S.; MacLane, S.
\paper Cohomology theory in abstract groups, II. Groups extensions 
with non-abelian kernel
\jour Ann. Math. (2)
\vol 48
\yr 1947
\pages 326--341
\endref

\ref
\key \cit0{4}
\by Giraud, Jean
\book Cohomologie non ab\'elienne
\bookinfo Grundlehren 179
\publ Springer-Verlag
\publaddr Berlin etc.
\yr 1971
\endref

\ref
\key \cit0{5}
\by Hochschild, G.
\paper Cohomology clases of finite type and finite dimensional 
kernels for Lie algebras 
\jour Am. J. Math.
\vol 76
\yr 1954
\pages 763-778
\endref

\ref
\key \cit0{6}
\by Hochschild, G. P.; Serre, J.-P.
\paper Cohomology of group extensions
\jour Trans. AMS 
\vol 74
\yr 1953
\pages 110--134
\endref

\ref
\key \cit0{7}
\by Hochschild, G. P.; Serre, J.-P.
\paper Cohomology of Lie algebras
\jour Ann. Math.
\vol 57
\yr 1953
\pages 591--603
\endref

\ref
\key \cit0{8}
\by Kac, V. G.
\paper Lie superalgebras
\jour Advances in Math.
\vol 26 
\yr 1997
\pages 8-96
\endref

\ref 
\key \cit0{9}
\by Kol\'a\v r, I.; Michor, Peter W.; Slov\'ak, J. 
\book Natural operations in differential geometry 
\publ Springer-Verlag
\publaddr Berlin Heidelberg New York
\yr 1993
\endref

\ref
\key \cit0{10}
\by Lecomte, Pierre
\paper Sur la suite exacte canonique asoci\'ee \`a un fibr\'e principal
\jour Bul. Soc. Math. France
\vol 13
\yr 1985
\pages 259--271
\endref
     
\ref
\key \cit0{11}
\by Lecomte, P.
\paper On some Sequence of graded Lie algebras associated to manifolds
\jour Ann. Global Analysis Geom. 
\vol 12
\yr 1994
\pages 183--192
\endref

\ref  
\key \cit0{12}  
\by Lecomte, Pierre; Michor, Peter W.; Schicketanz, Hubert 
\paper The multigraded Nijen\-huis-\-Richardson Algebra, its  
     universal property and application  
\jour J. Pure Applied Algebra 
\vol 77 
\yr 1992 
\pages 87--102 
\endref 

\ref
\key \cit0{13}
\by Lecomte, P.; Roger, C.
\paper Sur les d\'eformations des alg\`ebres de courants de type r\'eductif
\jour C. R. Acad. Sci. Paris, I
\vol 303
\yr 1986
\pages 807--810 
\endref

\ref
\key \cit0{14}
\by Mackenzie, K.
\book Lie groupoids and Lie algebroids in diferential geometry
\bookinfo London Mathematical Society Lecture Note Series, 124
\publ Cambridge University Press
\yr 1987
\endref

\ref
\key \cit0{15}
\by Mori, Mitsuya
\paper On the thre-dimensional cohomology group of Lie algebras
\jour J. Math. Soc. Japan
\vol 5
\yr 1953
\pages 171-183
\endref

\ref
\key \cit0{16}
\by Scheunert, M.
\book The theory of Lie superalgebras. An introduction
\bookinfo Lecture Notes in Mathematics. 716
\publ Springer-Verlag
\publaddr Berlin
\yr 1979
\endref

\ref
\key \cit0{17}
\by Serre, J.-P.
\paper Cohomologie des groupes discrets
\jour Ann. of Math. Studies
\vol 70
\yr 1971
\pages 77--169
\finalinfo Princeton University Press
\endref

\ref
\key \cit0{18}
\by Serre, J.-P.
\paper Cohomologie des groupes discrets
\jour S\'eminaire Bourbaki
\vol 399
\yr 1970/71
\endref

\ref
\key \cit0{19}
\by Shukla, U.
\paper A cohomology for Lie algebras
\jour J. Math. Soc. Japan
\vol 18
\yr 1966
\pages 275-289
\endref

\endRefs
\enddocument